\documentclass[11pt]{article}
\usepackage{amsmath}
\usepackage{amsthm}
\usepackage{mathtools}
\usepackage{amssymb}

\usepackage{enumitem}
\usepackage{hyperref}
\hypersetup{pdfcreator=XeLaTeX with hyperref}

\usepackage{geometry}
\usepackage{setspace}
\title{\textbf{{\normalsize PARTIALLY SMOOTH UNIVERSAL TAYLOR SERIES ON PRODUCTS OF SIMPLY CONNECTED DOMAINS}}} 
\author{Giorgos Kotsovolis \\ \\ Department of Mathematics \\ National and Kapodistrian University of Athens}
\date{}

\begin{document}
\maketitle 

\begin{abstract}
\noindent
Using a recent Mergelyan type theorem, we show the existence of universal Taylor series on products of planar simply connected domains $\Omega_i$ that extend continuously on $\prod\limits_{i=1}^{d}(\Omega_i \cup S_i)$, where $S_i$ are subsets of $\partial\,\Omega_i$, open in the relative topology. The universal approximation occurs on every product of compact sets $K_i$ such that $C-K_i$ are connected and for some $i_0$ it holds $K_{i_0}\cap (\Omega_{i_0}\cup \overline{S_{i_0}})=\varnothing$. Furthermore, we introduce some topological properties of universal Taylor series that lead to the voidance of some families of functions. \\ 
\end{abstract} 
\noindent
\textit{AMS classification numbers:} 30K05, 32A05, 32A17, 32A30 \\ \\
\textit{Keywords and phrases:} Universal Taylor series, Mergelyan's theorem, product of planar sets, Baire's theorem, abstract theory of universal series, generic property. 

\section{Introduction}
In one complex variable, Mergelyan's Theorem [17] has allowed the study of universal Taylor series [2], [6], [7], [10], [11], [12], [15], [16]. However, the study of universal Taylor series in several complex variables is underdeveloped [3], [4]. Recently, [5] gave birth to a new Mergelyan's type theorem. In particular, if $K=\prod\limits_{i=1}^{d}K_i$, where $K_i$ are planar compact sets with connected complements and $f: K\rightarrow C$ is a continuous function such that $f\circ \Phi$ is holomorphic on the disk $D\subset\mathbb{C}$, for every injective holomorphic mapping $\Phi: D\rightarrow K$, then $f$ is uniformly approximable on $K$ by polynomials. This result was recently used in [9], where it is proved that there exist holomorphic functions on $\prod\limits_{i=1}^{d}\Omega_i$, where $\Omega_i$ are simply connected domains of $C$, that behave universally on all products of planar compact sets, with connected complements, disjoint from $\prod\limits_{i=1}^{d}\Omega_i$, and also that there exist smooth functions on $\prod\limits_{i=1}^{d}\overline{\Omega_i}$, where $\Omega_i$ are simply connected domains of $C$, such that $C-\overline{\Omega_i}$ is connected, which behave universally on all products of planar compact sets, with connected complements, disjoint from $\prod\limits_{i=1}^{d}\overline{\Omega_i}$. In the present paper, these results are improved by proving the existence of holomorphic functions on $\prod\limits_{i=1}^{d}\Omega_i$, where $\Omega_i$ are planar simply connected domains that extend smoothly on $\prod\limits_{i=1}^{d}(\Omega_i \cup S_i)$,  where $S_i$ are subsets of $\partial\,\Omega_i$, open in the relative topology of $\partial S_i$ and $C-(\Omega_i\cup S_i)$ is connected. The universal approximation is achieved on all products of planar compact sets, with connected complements, disjoint from $\prod\limits_{i=1}^{d} (\Omega_i \cup \overline{S_i})$. Partially smooth universal Taylor series are not so developed even in one complex variable ([8], [19]). Our method cannot give the result for all compact sets $\prod\limits_{i=1}^{d} K_i$, $K_i \subset C$, $C-K_i$ connected which are disjoint from $\prod\limits_{i=1}^{d}(\Omega_i\cup S_i)$. In fact, we prove that there is no absorbing sequence of the previous sets. 

\section{Preliminaries} 
We will now present some abstract theory context from [1], [14], which will be of paramount importance for the main result of the paper. Let $(X_{\tau}), \,t\geq 1,$ indicate a sequence of metrizable  topological vector spaces over the field $\mathbb{K} \,(\mathbb{K}=\mathbb{R}$ or $\mathbb{K}=\mathbb{C})$, equipped with translation-invariant metrics $\rho_{\tau}$. Let $E$ denote a complete metrizable topological vector space, with topology induced by a translation-invariant metric $d$. Let, also, $L_m, \, m\geq1,$ be an increasing sequence of non-empty compact subsets of the metric space $L=\bigcup\limits_{m\geq 1} L_m$ and $\xi_0$ a distinguished element of $L_1$. Suppose that for any $n\in \mathbb{N}$, there exist continuous maps 
$$e_n:\, L\rightarrow E, \quad \chi_{\tau, n}:\, L\rightarrow X_{\tau} \quad \text{and} \quad \phi_n: L\times E \rightarrow \mathbb{K}.$$ 
Let $G$ denote the set of sequences in $\mathbb{K}^N$ with finitely many non-zero elements and let us suppose the following hold: 

\begin{enumerate}[label=(\roman*)]
\item The set $\{g_a: \, a\in G\}$ is dense in $E$, where $g_a=\sum\limits_{j=0}^{\infty}a_j e_j(\xi_0)$ and $a=(a_j)_{j=0}^{\infty}$. 
\item For every $a\in G$ and $\xi \in L$, the sets $\{n: \, \phi_n(\xi, g_a)\neq 0\}$ are finite and uniformly bounded with respect to $\xi\in L_m$ for any $m\geq 1$. 
\item For every $a\in G$ and $\xi \in L$: 
$$\sum\limits_{j=0}^{\infty}\phi_j(\xi, g_a)\,e_j(\xi) = \sum\limits_{j=0}^{\infty}a_j\, e_j(\xi_0)$$
\item For every $a\in G$, $\xi\in L$ and $\tau\geq 1$: 
$$\sum\limits_{j=0}^{\infty}\phi_j(\xi, g_a)\,\chi_{\tau, j}(\xi)=\sum\limits_{j=0}^{\infty}a_j\,\chi_{\tau, j}(\xi_0)$$ \end{enumerate}

\noindent
\textbf{Definition 2.1} Let $\mu$ be an infinite subset of $\{0, 1, 2, \dots\}$. Under the assumptions made above, an element $f\in E$ belongs to the class $U_{E, L}^{\mu}$, if, for every $\tau\geq 1$ and every $x\in X_{\tau}$, there exists a sequence $\lambda=(\lambda_n)\subset \mu $, such that, for every $m\geq1$ we have 
\begin{enumerate}[label=(\roman*)]
\item For every $\tau\geq1$, $$\sup\limits_{\xi\in L_m}\, \rho_{\tau} \left( \sum\limits_{j=0}^{\lambda_n}\phi_j (\xi, f)\,\chi_{\tau, j}(\xi), x \right) \rightarrow 0 \quad \text{as} \,\, u\rightarrow \infty$$ 
\item $$\sup\limits_{\xi\in L_m}\,d\left( \sum\limits_{j=0}^{\lambda_n}\phi_j (\xi, f)\,e_j(\xi), f \right) \rightarrow 0 \quad \text{as} \,\, u\rightarrow \infty.$$ 
\end{enumerate}
In the case $\mu=\{0, 1, 2, \dots\}$ we write $U_{E, L}$ instead of $U_{E, L}^{\mu}$. \\ \\ 
\textbf{Remark 2.2} In Definition 2.1 it is equivalent to ask that the sequence $(\lambda_n)\subset \mu $ is strictly increasing. [18] \\ \\
\textbf{Theorem 2.3} [1] \textit{Under the assumptions made above, the following are equivalent. 
\begin{enumerate}
\item For every $\tau\geq1$, every $x \in X_{\tau}$ and every $\varepsilon>0$, there exist $n\in\mathbb{N}$ and $a_0, a_1, a_2, \dots, a_n \in \mathbb{K}$ such that 
$$\rho_{\tau}\left( \sum\limits_{j=0}^{n}a_j\,\chi_{\tau, j}(\xi_0)\, ,\, x \right)<\varepsilon \quad \text{and} \quad d\left( \sum\limits_{j=0}^{n}e_j(\xi_0)\,a_j\, ,\, 0\right)<\varepsilon$$
\item For any increasing sequence $\mu$ of positive integers, the set $U_{E, L}^{\mu}$ is $G_{\delta}$ dense in $E$ and contains a dense vecotr subspace except zero. 
\end{enumerate} }
We now proceed in defining some classes of functions that will be used widely in the present paper. We will also present some new results from [5], in forms of lemmas, that will prove useful later on. \\ \\ 
\textbf{Definition 2.4} Recalling that a mapping from a planar domain to $C^d$ is holomorphic, if each coordinate is a complex valued holomorphic function, we define the algebra $A_D(K)$ as the set of all functions $f:\, K\rightarrow C$, where $K$ is compact in $C^d$, that are continuous on $K$ and such that, for every open disc $D\subset C$ and every holomorphic mapping $\phi: D\rightarrow K$, the composition $f\circ \phi: D\rightarrow C$ is holomorphic. \\ \\ 
It is easy to see that $A_D(K)$ becomes a metrizable topological vector space, with the supremum norm. \\ \\  
\textbf{Lemma 2.5} \textit{Let $K=\prod\limits_{i=1}^{d}K_i$, where $K_i\subset C$ are compact sets, with $C-K_i$ connected. If $f\in A_D(K)$ and $\varepsilon>0$, there exists a polynomial $P$ such that $\sup\limits_{z\in K}\vert f(z)-P(z)\vert <\varepsilon$.} \\ \\ 
Lemma 2.5 is part of a result in [5], Remark 4.9 (3). From the same paper, we also get the following. \\ \\
\textbf{Definition 2.6} Let $K$ be a compact subset of $C^d$. We define the algebra $O(K)$ as the set of all functions $f: K\rightarrow C$, such that there exists an open set $U\supseteq K$, so that $f\in H(U)$. \\ \\ 
It is easy to see that, if we enduce the algebra $O(K)$ with the denumerable family of seminorms $\Big\{\sup\limits_{z\in K}\vert D_a f\vert,\,\, D_a \,\,$ varies in the denumerable set of all differential operators of mixed partial derivatives in $z=(z_1, z_2, \dots, z_d) \Big\}$, then $O(K)$ becomes a Fr\'echet space. \\ \\ 
\textbf{Lemma 2.7} \textit{Let $K$ be as in Lemma 2.5. If $f\in O(K)$, $I$ is a finite subset of $N^d$ and $\varepsilon>0$, there exists a polynomial $Q$ such that 
$$\sup\limits_{z\in K} \bigg\vert \dfrac{\partial^{a_1+\dots+a_d}}{\partial_{z_1}^{a_1}\,\partial_{z_2}^{a_2}\dots \partial_{z_d}^{a_d}}(f-Q)(z)\bigg\vert< \varepsilon$$ for every $a=(a_1, a_2, \dots, a_d)\in I$. } \\ \\ 
\textbf{Remark 2.8} Namely, Lemmas 2.5 and 2.7 inform us that, if $K=\prod\limits_{i=1}^{d} K_i$, where $K_i$ are compact sets of $C$ and $C-K_i$ is connected, then the set of polynomials is dense in both $A_D(K)$ and $O(K)$. \\ \\ 
The following definition will be the one that will most concern us in the current paper. \\ \\ 
\textbf{Definition 2.9} Let $\Omega_i \subsetneq C, \, 1\leq i\leq d,$ be simply connected domains of $C$ and $\Omega=\prod\limits_{i=1}^{d}\Omega_i$. Let also $S_i$ be a subset of $\partial\,\Omega_i$, for each $i\in\{1, 2, \dots, d\}$. If $S=\prod\limits_{i=1}^{d} S_i$, we define the class $A^{\infty}(\Omega, S)$ of all the functions $f: \Omega \rightarrow C$ holomorphic on $\Omega$, such that for every differential operator $D$ of mixed partial derivatives $Df$ is continuously extendable on $\prod\limits_{i=1}^{d}(S_i\cup \Omega_i)$. \\ \\ 
\textbf{Theorem 2.10} \textit{Along with the assumptions of Definition 2.9, we also assume that $C-(S_i\cup \Omega_i)$ are connected, for $1\leq i \leq d$. Then, there exists a sequence $L_n$ in $\prod\limits_{i=1}^{d}(S_i\cup\Omega_i)$, $n\geq 1,$ such that 
\begin{enumerate}[label=(\roman*)]
\item $L_n$ are products of planar compact sets with connected complements. 
\item $L_n$ is increasing with respect to $n$. 
\item For every compact set $L$ in $\prod\limits_{i=1}^{d} (\Omega_i\cup S_i)$, there exists $n\geq1$, such that $L\subseteq L_n$.
\end{enumerate}
Taking the above facts into consideration, if $A^{\infty}(\Omega, S)$ is enduced with the denumerable family of seminorms $\Big\{\sup\limits_{z\in L_n}\vert D_a (f)(z)\vert,\,\, n\geq 1\,\, and \,\,D_a \,\,$ varies in the denumerable family of differential operators of mixed partial derivatives in $z=(z_1, z_2, \dots, z_d) \Big\}$, it becomes a Fr\'echet space.} \\ \\ 
\textit{Proof.} Let $L$ be a compact set in $\prod\limits_{i=1}^{d} (\Omega_i\cup S_i)$. Without loss of generality, since $\Omega$ is a product of planar sets, we may assume $L=\prod\limits_{i=1}^{d} L'_i$, where $L'_i$ is a compact set of $\Omega_i \cup S_i$. The set $\partial\,\Omega_i \setminus S_i$ is closed and, thus, there exists an $n$ such that $L'_i\subseteq \Big\{ z\in C: \text{dist}(z\, ,\, \partial\,\Omega_i\setminus S_i)\geq \dfrac{1}{n} \Big\}$. We can also demand that $$L'_i\subseteq \overline{\Omega_i}\, \cap\, \overline{B}(0, n)\,\cap\, \Big\{ z\in C: \text{dist}(z\, ,\, \partial\,\Omega_i\setminus S_i)\geq \dfrac{1}{n} \Big\}.$$ Let $L_{i, n}=\overline{\Omega_i}\, \cap\, \overline{B}(0, n)\,\cap\, \Big\{ z\in C: \text{dist}(z, \partial\,\Omega_i\setminus S_i)\geq \dfrac{1}{n} \Big\}$. We proved that for every compact set $L'_i$ in $\Omega_i\cup S_i$ there exists $n\in\mathbb{N}$ such that $L'_i\subset L_{i, n}$. It is trivial to see that $L_{i, n}$ is a compact set of $\Omega_i\cup S_i$. Now, we will prove that
$$(L_{i, n})^c = \Big(\overline{\Omega_i}\Big)^c \, \cup\, \Big(\overline{B}(0, n)\Big)^c\,\cup\, \Big\{ z\in C: \text{dist}(z, \partial\,\Omega_i\setminus S_i)< \dfrac{1}{n} \Big\}$$ is a connected set. It is sufficient to show that $\Big(\overline{\Omega_i}\Big)^c \, \cup\, \Big\{ z\in C: \text{dist}(z, \partial\,\Omega_i\setminus S_i)< \dfrac{1}{n} \Big\}$ is connected and that $\Big(\overline{\Omega_i}\Big)^c \, \cup\, \Big\{ z\in C: \text{dist}(z, \partial\,\Omega_i\setminus S_i)< \dfrac{1}{n} \Big\}$ is not disjoint from $\Big(\overline{B}(0, n)\Big)^c$. For the former, notice that $\Big(\overline{\Omega_i}\Big)^c \, \cup\, \Big\{ z\in C: \text{dist}(z, \partial\,\Omega_i\setminus S_i)< \dfrac{1}{n} \Big\}$ is equal to the set $$(C- (\Omega_i\cup S_i))\, \cup\, \Big\{ z\in C: \text{dist}(z, \partial\,\Omega_i\setminus S_i)< \dfrac{1}{n} \Big\},$$ which can be written down as 
$$(C- (\Omega_i\cup S_i))\, \cup\, \Bigg( \bigcup\limits_{z\in\partial(\Omega_i)\setminus S_i} B\left(z, \dfrac{1}{n} \right) \Bigg) $$ which is connected, since $B\left(z, \dfrac{1}{n} \right)$ is connected and $$z\in(C- (\Omega_i\cup S_i))\, \cap\, B\left(z, \dfrac{1}{n} \right) \neq\varnothing $$ for each $z\in\partial(\Omega_i)\setminus S_i$. For the latter, note that if the set $\Big(\overline{\Omega_i}\Big)^c \, \cup\, \Big\{ z\in C: \text{dist}(z, \partial\,\Omega_i\setminus S_i)< \dfrac{1}{n} \Big\}$ is disjoint from $\Big(\overline{B}(0, n)\Big)^c$, then so is the set 
$C- (\Omega_i\cup S_i)$. But then $\Big(\overline{B}(0, n)\Big)^c\subseteq \Omega_i\cup S_i$ and, thus, 
$$\Big(\overline{B}(0, n)\Big)^c = \left(\Big(\overline{B}(0, n)\Big)^c\right)^{\circ}\subseteq (\Omega_i\cup S_i)^{\circ}\subseteq \Omega_i. $$ Since $\Omega_i$ is simply connected, it follows that $\Omega_i=C$, which is absurd. Notice now that $L_n = \prod\limits_{i=1}^{d} L_{i, n}$ is the desired sequence. \qed \\ \\ 
\textbf{Definition 2.11} We define $X^{\infty}(\Omega, S)$ the closure of the set of polynomials in the Fr\'echet space $A^{\infty}(\Omega, S)$. 

\section{Partially smooth universal Taylor series}
Let $\Omega\subseteq C^d$ be an open set, $\zeta=(\zeta_1, \zeta_2, \dots, \zeta_d)$ be a point of $\Omega$ and $N_j$, $j=1, 2, \dots$ be an enumeration of $N^d$, where $N=\{0, 1, 2, \dots\}$, and denote the monomials $$(z-\zeta)^a = (z_1-\zeta_1)^{a_1} (z_2-\zeta_2)^{a_2}\cdots (z_d-\zeta_d)^{a_d}$$ where $a=(a_1, a_2, \dots, a_d)$. If $f$ is a holomorphic function on $\Omega$, we denote by $a(f, \zeta)$ the coefficient of $(z-\zeta)^a$ in the Taylor expansion of $f$ with center $\zeta$, that is $$a(f, \zeta)=\dfrac{1}{a_1!\,a_2!\cdots a_d!} \dfrac{\partial^a}{\partial_{z_1}^{a_1}\,\partial_{z_2}^{a_2}\dots \partial_{z_d}^{a_d}}\,f(\zeta).$$ Let $S_N(f, \zeta)(z)=\sum\limits_{j=0}^{N} a_{N_j}(f, \zeta)\,(z-\zeta)^{N_j}$. \\ \\ 
\textbf{Theorem 3.1} \textit{Let $\Omega_i$, $i=1, 2, \dots, d$, be simply connected domains and $\Omega=\prod\limits_{i=1}^{d}\Omega_i$. Let also $S_i \subset \partial\,\Omega_i$ be open sets in the topology of $\partial\,\Omega_i$, such that $C-(\Omega_i \cup S_i)$ is connected for $1\leq i\leq d$, and $S=\prod\limits_{i=1}^{d}S_i$. Fix $\zeta_0\in \Omega$ and the enumeration $N_j$ as above. Let $\mu$ be n infinite subset of $N$. There exists a function $A^{\infty}(\Omega, S)$, such that for every product of planar compact sets $K=\prod\limits_{i=1}^{d}K_i$, where $K_i$ have connected complements and $K$ is disjoint from $\prod\limits_{i=1}^{d}(\Omega_i\cup \overline{S_i})$, and every $h\in A_D(K)$, there exists a strictly increasing sequence $\lambda_n\in \mu$, $n=1, 2, \dots$, such that $S_{\lambda_n}(f, \zeta^{\circ})(z)\rightarrow h(z)$ uniformly on $\prod\limits_{i=1}^{d}K_i$ and $S_{\lambda_n}(f, \zeta^{\circ})(z)\rightarrow f(z)$ in the topology of $A^{\infty}(\Omega, S)$, as $n\rightarrow +\infty$. Furthermore, the set of such $f\in X^{\infty}(\Omega, S)$ is a dense $G_{\delta}$ set and contains a dense vector subspace except zero.} \\ \\ 
\textit{Proof.} According to the Remark 2.2, it is not necessary to prove that the sequence $\lambda_n$ is increasing; it suffices that $\lambda_n\in\mu$ for all $n$. \par It is known [12], [16], that when $\Omega_i$ is a (simply) connected domain in $C$, there can be found compact sets $K_{i, 1}, K_{i, 2}, \dots$ in $C-\Omega_i$, with $C-K_{i, j}$ connected, such that for every $K$ compact set in $C-\Omega_i$, with $C-K$ connected, there exists $j$ such that $K\subseteq K_{i, j}$. Let 
$$\widetilde{K}_i = \bigg\{K_{i,j}\cap\Big\{z: \text{dist}(z, \overline{S_i})\geq\dfrac{1}{m}\Big\},\,\,\text{for}\,\,j\in N\,\, \text{and}\,\, m\in N\bigg\}.$$
Let $\widetilde{K}_{i,j}$ be an enumeration of the set $\widetilde{K}_i$. One can also verify that 
\begin{enumerate}[label=(\roman*)]
\item $\widetilde{K}_{i,j}$ is compact, with $C-\widetilde{K}_{i,j}$ connected for $1\leq i \leq d$ and $j\in\{1, 2, 3, \dots\}$. 
\item If $K$ is a compact set in $C-(\Omega_i \cup \overline{S_i})$, with $C-K$ connected, there exists $j$ such that $K\subseteq K_{i,j}$.
\end{enumerate}
Let $K_{\tau}$, $\tau=1, 2, \dots,$ be an enumeration of all $\prod\limits_{i=1}^{d} T_i$, where $T_{i_0}=\widetilde{K}_{i_0, j_{i_0}}$ for some $i_0$ and $j_{i_0}$ and $T_i=\overline{B}(0, s_i)$ for all $i\neq i_0$ and $s_i$ varies in the set of natural numbers. We notice that if we have the theorem for compact sets $K$ of the form $K_{\tau}$, we also have it for every $K=\prod\limits_{i=1}^{d}V_i$, where $V_i$ are compact sets with connected complements and $K\cap \left( \prod\limits_{i=1}^{d}(\Omega_i\cup\overline{S_i}) \right)=\varnothing$. Indeed, let us assume that the function $f$ satisfies the theorem for compact sets $K$ of the form $K_{\tau}$ and $h\in A_D\left(\prod\limits_{i=1}^{d}V_i\right)$. Let also $\varepsilon>0$. According to Lemma 2.5, there exists a polynomial $P$ such that $$\sup\limits_{z\,\in\prod\limits_{i=1}^{d}V_i}\vert P(z)-h(z)\vert < \dfrac{\varepsilon}{2}.$$ By definition of the sets $K_{\tau}$, there exists $\tau_0\in N$, such that $\prod\limits_{i=1}^{d}V_i\subset K_{\tau_0}$. The theorem is valid for the compact set $K_{\tau_0}$ and thus there exists $\lambda\in\mu$, such that 
$$\sup\limits_{z\, \in K_{\tau_0}} \big\vert S_{\lambda}(f, \zeta^{\circ})(z) - P(z) \big\vert <\dfrac{\varepsilon}{2} \quad\text{and}\quad d\left(S_{\lambda}(f, \zeta^{\circ})\, , \, f\right)<\varepsilon$$ where $d$ is the metric of the Fr\'echet space $A^{\infty}(\Omega, S)$. The result follows by the triangular inequality. \\ Now, we proceed to show the result for the compact sets $K_{\tau}$. Following the notation of Theorem 2.3, we set $X_{\tau}=A_D(K_{\tau})$, $E=X^{\infty}(\Omega, S)$, $L_m=\{\zeta_0\}$ for every $m\in N$ and $\xi_0=\zeta_0$. Let, also, $$e_n=\chi_{\tau, n}(\xi)=\left(z\rightarrow (z-\xi)^{N_n}\right) \quad \text{and} \quad \phi_n(\xi, f)=a_{N_n}(f, \xi).$$
The conditions (i), (ii), (iii) and (iv) of Section 2 can easily be verified. According, now, to Theorem 2.3, it suffices to show that, for every $t\geq 1$, $g\in A_D(K_{\tau})$, every $\varepsilon>0$, every $\widetilde{L}\subseteq \prod\limits_{i=1}^{d}(\Omega_i\cup S_i)$ compact and every finite set $F$ of mixed partial derivatives in $z=(z_1, z_2, \dots, z_d)$, there exists a polynomial $P$ such that $$\sup\limits_{z\,\in K_{\tau}}\vert P(z) - g(z)\vert <\varepsilon \quad \text{and} \quad \sup\limits_{z\,\in \widetilde{L}}\vert D_{\ell}(P)(z)\vert <\varepsilon$$ for every $D_{\ell}\in F$. Since $g\in A_D(K_{\tau})$, there exists by Lemma 2.5, a polynomial $\tilde{g}$ such that $$\sup\limits_{z\,\in K_{\tau}}\vert g(z) - \tilde{g}(z)\vert <\dfrac{\varepsilon}{2}.$$  
Furthermore, since $\prod\limits_{i=1}^{d}(\Omega_i\cup S_i)$ is a product of planar sets, we can, without loss of generality, assume that $\widetilde{L}=\prod\limits_{i=1}^{d}\widetilde{L}_i$, where $\widetilde{L}_i$ are compact sets in $\Omega_i \cup S_i$ and $C-\widetilde{L}_i$ are connected. Denote, as well, $K_{\tau}=\prod\limits_{i=1}^{d}K_{i, \tau}$, where $K_{i, \tau}$ are all balls except for $i=i_0$. \\ 
For $i\neq i_0$ there exist balls $\overline{B_i}$, such that $K_{i, \tau}\cup \widetilde{L}_i\subseteq \overline{B_i}$. We define a function $h:\, \prod\limits_{i=1}^{d} K_{i, \tau}\cup \widetilde{L}_i\rightarrow C$ as follows: $h(z)=\tilde{g}(z)$ for $z\,\in \prod\limits_{i=1}^{d}Q_i$, where $Q_{i_0}=K_{i_0, \tau}$ and $Q_i=\overline{B_i}$ for $i\neq i_0$, and $h(z)=0$ for $z\,\in \prod\limits_{i=1}^{d}Q'_i$, where $Q'_i=\overline{B_i}$ for $i\neq i_0$ and $Q'_{i_0}=\widetilde{L}_{i_0}$.  \\
Notice that the function can be extended on an open neighborhood of $\prod\limits_{i=1}^{d}(Q_i\cup Q'_i)$ and thus, by Lemma 2.7 applied for the function $h$ and $\dfrac{\varepsilon}{2}$, the wanted polynomial exists. The proof is complete. \qed \\ \\ 
\textbf{Remark 3.2} The last part of the above proof shows that the result of Theorem 3.1 holds for any denumerable family of compact sets $K=\prod\limits_{i=1}^{d}K_i$, $K_i\subset C$, $C-K_i$ connected, which are disjoint from $\prod\limits_{i=1}^{d}(\Omega_i\cup S_i)$. In particular, this holds for any singleton in $\left( \prod\limits_{i=1}^{d}(\Omega_i\cup S_i)\right)^c$. \\ \\ 
\textbf{Remark 3.3} Our method has succeeded in proving universal behaviour on the products of compact sets with connected complements, which are disjoint from $\prod\limits_{i=1}^{d}(\Omega_i\cup \overline{S_i})$. The question arises whether universal behaviour can occur on the products of compact sets with connected complements in the complement of $\prod\limits_{i=1}^{d}(\Omega_i\cup S_i)$, or not. Namely, can we partition $C^d$ in a "smooth" and a "universal" set with respect to a function $f$? Of course, if $S_i$ is clopen in the topology of $\partial\,\Omega_i$, for every $i=1, 2, \dots, d$, the answer is positive by Theorem 3.1, because $S_i=\overline{S_i}$. Take for instance $\Omega$ to be a strip and $S$ be one of the two connected components of $\partial\,\Omega$. The following Lemma shows that our method fails in general. \\ \\ 
\textbf{Lemma 3.4} \textit{Let $M$ be a set in $C$, such that there exists a sequence $K_i$, $i\geq 1$, of compact subsets of $M$, such that for every compact set $K$ in $M$, with connected complement $C-M$, there exists $i_0\in N$ so that $K\subseteq K_{i_0}$. Then, $M\cap \partial M$ is open in the topology of $\partial M$. } \\ \\ 
\textit{Proof.} Suppose the contrary. Then, there exists a sequence $z_i, \,i\geq 1$ in $\partial M \setminus (\partial M \cap M)$, such that $z_i$ converges to $z$, a point of $\partial M \cap M$. Since $z_i\,\in K_i$ and $K_i$ is compact, for every $i\in N$, we can find a point $t_i$ of $M\setminus K_i$ such that dist$(t_i, z_i)<\dfrac{1}{i}$. Note that $t_i$ converges on $z$. Let $K=\{t_i, \,\,i\in N\}\cup\{z\}$. Then, $K$ is compact with $C-K$ connected, but it is not contained in any $K_i$, $i\in N$, which is a contradiction. Thus, $\Omega \cap \partial\,\Omega$ must be open in the topology of $\partial\,\Omega$. \qed \\ \\ 
\textbf{Remark 3.5} Notice that if we try to show that the universal behaviour occurs on the products of compact sets with connected complements, in the complement of $\prod\limits_{i=1}^{d}(\Omega_i\cup S_i)$, it is essential for our method to find a denumerable family of compact sets in $C^d-\prod\limits_{i=1}^{d}(\Omega_i\cup S_i)$, $K_i\geq 1$, such that $K_i$ are products of compact sets with connected complements and for every $K$ in $C^d-\prod\limits_{i=1}^{d}(\Omega_i\cup S_i)$, product of compact sets with connected complements, there exists $i_0\in N$, so that $K\subseteq K_{i_0}$. Let $z_i, \, i=2, \dots, d$, be points in $\Omega_i\cup S_i$. If $K'$ is a compact set in $C-(\Omega_1\cup S_1)$, with connected complement, there can be found $i_0\in N$ such that $$K'\times \{z_2\}\times \dots \times \{z_d\} \subseteq K_{i_0}.$$ 
Denote $K_{i_0}=\prod\limits_{j=1}^{d} K_{i_0, j}$ and notice that $z_j\in K_{i_0, j}$ for $j\geq 2$. Thus, since $K_{i_0}\subseteq C^d- \prod\limits_{i=1}^{d}(\Omega_i\cup S_i)$ and $K_{i_0, 1}\subseteq C-(\Omega_1\cup S_1)$, therefore, if $\widetilde{K}=\big\{ K_{i,1},\,\, \text{for} \,\, K_{i,1}\subseteq C-(\Omega_i\cup S_i) \big\}$ and $K_i$ is an enumeration of of the elements of the set $\widetilde{K}$, we deduce that for every compact set $K$ in $C-(\Omega_1\cup S_1)$, with connected complement, there exists $i_0\in N$ so that $K\subseteq \widetilde{K}_{i_0}$. By the previous lemma, $S_1$ is closed in the topology of $\partial\,\Omega_1$. The same, of course, holds for all $i\in \{1, 2, \dots, d\}$. Thus, we conclude that $S_i$ is clopen in $\partial\,\Omega_i$, for all $i\in \{1, 2, \dots, d\}$. \\ \\
\textbf{Remark 3.6} In the extreme case where $S_i=\varnothing$ for all $i$, or $S_i=\partial\,\Omega$ for all $i$, we get the results of [9], regarding holomorphic and smooth, respectively, universal Taylor series. \\ \\ 
\textbf{Theorem 3.7} \textit{Under the assumptions and notation of Theorem 3.1, there exists a function $f\in A^{\infty}(\Omega, S)$, such that for every planar compact sets $K_i$ with $C-K_i$ connected and $\left(\prod\limits_{i=1}^{d}K_i \right)\cap \left(\prod\limits_{i=1}^{d}(\Omega_i\cup \overline{S_i})\right)=\varnothing$, and every $h\in A_D\left(\prod\limits_{i=1}^{d}K_i\right)$, there exists a strictly increasing sequence $\lambda_n\in \mu$ such that for every compact set $\widetilde{L}\subseteq\prod\limits_{i=1}^{d}(\Omega\cup S_i)$ we have $$\sup\limits_{z\,\in \prod\limits_{i=1}^{d}(K_i)\, ,\, \zeta\in\widetilde{L}} \vert S\lambda_n(f, \zeta)-h(z) \vert\rightarrow 0 \quad \text{and} \quad \sup\limits_{\zeta\in \widetilde{L}} d\left(S\lambda_n(f, \zeta)\, , \, f\right)\rightarrow 0$$ as $n\rightarrow+\infty$, where $d$ is the metric of the Fr\'echet space $A^{\infty}(\Omega, S)$. Furthermore, the set of such functions $f\in X^{\infty}(\Omega, S)$ is $G_{\delta}$-dense and contains a dense vector subspace except zero.} \\ \\ 
\textit{Proof.} The proof is analogous to the proof of Theorem 3.1. The only difference is that when we apply Theorem 2.3, we do not set $L_m=L=\{\zeta_0\}$, but instead we set $L=\prod\limits_{i=1}^{d}(\Omega_i\cup S_i)$ and $L_m, \,m=1, 2, \dots,$ the family of compact subsets of $\prod\limits_{i=1}^{d}(\Omega_i\cup S_i)$, discussed in Theorem 2.10, and $\zeta_0$ an arbitrary element of $L_1$. \qed \\ \\ 
We now proceed by showing similar theorems with the ones just proved, but this time the universal approximation occurs on functions of the space $O(K)$, instead of the space $A_D(K)$. \\ \\ 
\textbf{Theorem 3.8} \textit{Under the assumptions and notation of Theorem 3.1, there exists a function $f\in A^{\infty}(\Omega, S)$, such that for every planar compact sets $K_i$ with $C-K_i$ connected and $\left(\prod\limits_{i=1}^{d}K_i \right)\cap \left(\prod\limits_{i=1}^{d}(\Omega_i\cup \overline{S_i})\right)=\varnothing$, and every $h\in O\left(\prod\limits_{i=1}^{d}K_i\right)$, there exists a strictly increasing sequence $\lambda_n\in \mu$ such that $S\lambda_n(f, \zeta^{\circ})(z)\rightarrow h(z)$ in the topology of $O\left(\prod\limits_{i=1}^{d}K_i\right)$ and $S\lambda_n(f, \zeta^{\circ})(z)\rightarrow f(z)$ in the topology of $A^{\infty}(\Omega, S)$, as $n\rightarrow+\infty$. Furthermore, the set of such functions $f\in X^{\infty}(\Omega, S)$ is a $G_{\delta}$-dense set and also contains a dense vector subspace except zero.} \\ \\ 
\textit{Proof.} According to Remark 2.2, it is sufficient to prove the theorem not for $\lambda_n$ strictly increasing, but for a sequence $\lambda_n\in\mu$ for all $n$. \\ 
Let $K_{\tau}$, $\tau=1, 2, \dots$, represent the sequence of compact sets of $C^d$ , discussed in the proof of Theorem 3.1. Notice, now, that if the theorem is valid for compact sets of the form $K_{\tau}$, $\tau=1, 2, \dots$, we also have it proved for every $K=\prod\limits_{i=1}^{d}T_i$ disjoint from $\prod\limits_{i=1}^{d}(\Omega_i\cup \overline{S_i})$, such that $C-T_i$ are connected and $T_i$ are compact. Indeed, let us suppose that the theorem is valid for compact sets $K$ of the form $K_{\tau}$. Let $h\in O\left( \prod\limits_{i=1}^{d}T_i\right)$, $\varepsilon>0$ and $F$ be a finite set of differential operators of mixed partial derivatives with respect to $z=(z_1, z_2, \dots, z_d)$. Using Lemma 2.7, we can find a polynomial $P$ such that 
$$\sup\limits_{z\,\in\prod\limits_{i=1}^{d}T_i}\vert DP(z)-Dh(z)\vert<\dfrac{\varepsilon}{2} \quad \text{for all}\,\, D\in F.$$
By definition of the sequence $K_{\tau}$, there exists $\tau_0\in N$, such that $\prod\limits_{i=1}^{d}T_i \subseteq K_{\tau_0}$. The theorem is valid for compact sets of that form and thus we can find $\lambda\in\mu$, so that 
$$\sup\limits_{z\,\in K_{\tau_0}}\big\vert DP(z)-D(S\lambda(f, \zeta^{\circ}))(z)\big\vert<\dfrac{\varepsilon}{2} \quad \text{for all}\,\, D\in F.$$ and 
$$d\left( S\lambda(f, \zeta^{\circ})\, , \, f\right)<\dfrac{\varepsilon}{2}$$ where $d$ is the metric of $A^{\infty}(\Omega, S)$. The result follows now by the triangular inequality.\par
We proceed now, using the above fact, in proving the theorem for compact sets $K$ of the form $K_{\tau}$. Following the notation of Theorem 2.3, we set  $X_{\tau}=O(K_{\tau})$, $E=X^{\infty}(\Omega, S)$, $L_m=\{\zeta_0\}$ for every $m\in N$ and $\xi_0=\zeta_0$.
The conditions (i), (ii), (iii) and (iv) of Section 2 can easily be verified. According to Theorem 2.3, it suffices to show that for every $\tau\geq1$, $g\in O(K_{\tau})$, every $\varepsilon>0$, every $\widetilde{L}$ compact subset of $\prod\limits_{i=1}^{d}(\Omega_i\cup S_i)$, and every finite set $F$ of differential operators of mixed partial derivatives with respect to $z=(z_1, z_2, \dots, z_d)$ there exists a polynomial $P$ such that 
$$\sup\limits_{z\,\in K_{\tau}}\big\vert DP(z)-Dg(z)\big\vert<\varepsilon \quad \text{and} \quad \sup\limits_{z\,\in \widetilde{L}}\vert DP(z)\vert <\varepsilon \quad \text{for all}\,\,D\in F.$$
Since $g\in O(K_{\tau})$, there exists by Lemma 2.7, a polynomial $\tilde{g}$, so that $$\sup\limits_{z\,\in K_{\tau}}\big\vert D\tilde{g}(z)-Dg(z)\big\vert<\dfrac{\varepsilon}{2} \quad \text{for all}\,\,D\in F.$$
Moreover, since $\prod\limits_{i=1}^{d}(\Omega_i\cup S_i)$ is a product of planar sets, we can assume, without loss of generality, that $\widetilde{L}=\prod\limits_{i=1}^{d}\widetilde{L}_i$, where $\widetilde{L}_i$ are compact subsets of $\Omega_i\cup S_i$, for each $1\leq i \leq d$, and $C-\widetilde{L}_i$ is connected (according to Theorem 2.10). Denote also $K_{\tau}=\prod\limits_{i=1}^{d}K_{i, \tau}$, where $K_{i, \tau}$ are all balls except for $i=i_0$. For $i\neq i_0$ there exist balls $\overline{B_i}$ such that $K_{i, \tau}\cup \widetilde{L}\subseteq \overline{B_i}$. \par 
We define a new function $u$ as follows: $u(z)=\tilde{g}(z)$ for $z\,\in \prod\limits_{i=1}^{d}Q_i$, where $Q_{i_0}=K_{i_0, \tau}$ and $Q_i=\overline{B_i}$ for $i\neq i_0$, and $u(z)=0$ for $z\,\in \prod\limits_{i=1}^{d}Q'_i$, where $Q'_{i_0}=\widetilde{L}_{i_0}$ and $Q'_{i_0}=\overline{B_i}$ for $i\neq i_0$. \par 
The function $u$ is well defined and it is easy to see that $u \in O\left(\prod\limits_{i=1}^{d}(Q_i\cup Q'_i)\right)$. The desired polynomial $P$ is now givem by Lemma 2.7 applied to the function $u$, $\dfrac{\varepsilon}{2}$ and set $F$ of differential operators. \qed \\ \\ 
\textbf{Theorem 3.9} \textit{Under the assumptions and notation of Theorem 3.1, there exists a function $f\in A^{\infty}(\Omega, S)$, such that for every compact set $K_i$ with connected complements, such that $K=\prod\limits_{i=1}^{d}K_i$ is disjoint from $\prod\limits_{i=1}^{d}(\Omega_i\cup \overline{S_i})$, and every $h\in O(K)$, there exists a strictly increasing sequence $\lambda_n\in \mu$, so that for every compact set $\widetilde{L}\subseteq \prod\limits_{i=1}^{d}(\Omega_i\cup S_i)$ we have $$\sup\limits_{\zeta\in \widetilde{L}}\big\vert \rho(S\lambda_n(f, \zeta)\, , \, h) \big\vert \rightarrow 0$$ where $\rho$ is the metric of the Fr\'echet space $O(K)$, and $$\sup\limits_{\zeta\in \widetilde{L}}\big\vert d(S\lambda_n(f, \zeta)\, , \, f) \big\vert \rightarrow 0$$ where $d$ is the metric of the Fr\'echet space $A^{\infty}(\Omega, S)$, as $n\rightarrow+\infty$. Furthermore, the set of such functions $f\in X^{\infty}(\Omega, S)$ is a $G_{\delta}$-dense set and contains, as well, a dense vector subspace except zero.} \\ \\ 
\textit{Proof.} The proof is similar to the proof of Theorem 3.8. The only difference occurs when Theorem 2.3 is applied. In detail, we do not set $L_m=L=\{\zeta_0\}$, but instead we set $L=\prod\limits_{i=1}^{d}(\Omega_i\cup S_i)$ and $L_m$, $m=1, 2, \dots$ the family of compact subsets of $\prod\limits_{i=1}^{d}(\Omega_i\cup S_i)$, discussed in Theorem 2.10, and $\zeta_0$ an arbitrary set of $L_1$. \qed \\ \\ 
\textbf{Remark 3.10} We notice that the previous theroems have been proved with a specific enumeration of the monomials in mind. The question whether there exist similar results simultaneously with respect to all enumerations of the monomials has been answered in the negative in [9], Theorem 6.1, for the case where $S_i=\varnothing$ for all $i=1, 2, \dots, d$. The argument can also be applied in our case. \\ \par
We will now present some applications of the theorems in this section. \\ \\ 
\textbf{Definition 3.11} If $F$ is a family of differential operators of mixed partial derivatives, we call $F$ \textit{gapless} if for every $\dfrac{\partial^{a_1+\dots+a_d}}{\partial_{z_1}^{a_1}\dots\partial_{z_d}^{a_d}}\in F$ and every sequence $b_i$ with $b_i\leq a_i$, it holds that $\dfrac{\partial^{b_1+\dots+b_d}}{\partial_{z_1}^{b_1}\dots\partial_{z_d}^{b_d}}$ is also contained in $F$. \\ \\
\textbf{Definition 3.12} (see also [13]) Let $\Omega_i\subsetneq C$, $1\leq i\leq d$, be simply connected domains in $C$ and $\Omega=\prod\limits_{i=1}^{d}\Omega_i$. Let also $S_i$ be a subset of $\partial\,\Omega_i$, for each $i\in\{1, 2, \dots, d\}$. Denote by $F$ a gapless family of differential operators of mixed partial derivatives. If $S=\prod\limits_{i=1}^{d}S_i$, we define the class $A^F(\Omega, S)$ of all the functions $f:\Omega\rightarrow S$, holomorphic on $\Omega$, such that for every $D\in F$, $Df$ is continuously extendable on $\prod\limits_{i=1}^{d}(\Omega_i\cup S_i)$. \\ \\ 
\textbf{Remark 3.13} Notice that, following the notification of Definition 3.12, if $S_i$ are open in the relative topology of $\partial\,\Omega_i$, then $A^F(\Omega, S)$ becomes a Fr\'echet space, with topology induced by the denumerable family of seminorms $\bigg\{ \sup\limits_{z\in L_n} \vert Df\vert, \,\,n\geq 1,\,D\in F \bigg\}$, where $L_n$, $n\geq 1$, is the sequence discussed in Theorem 2.10. \\ \\ 
\textbf{Definition 3.14} With the assumption of Remark 3.13, we define $X^F(\Omega, S)$ the closure of the set of the polynomials in the Fr\'echet $A^F(\Omega, S)$. \\ \\ 
\textbf{Lemma 3.15} [1] \textit{Under the assumptions and the notation of Theorem 2.3, the following are equivalent. 
\begin{enumerate}
\item For any increasing sequence $\mu$ of positive integers, the set $U_{E, L}^{\mu}$ is not empty. 
\item For any increasing sequence $\mu$ of positive integers, the set $U_{E, L}^{\mu}$ is $G_{\delta}$ dense in $E$ and contains a dense vector subspace except zero.
\end{enumerate} } 
\noindent
\textbf{Remark 3.16} Notice that the theorems of this section hold, also, if we substitute $A^{\infty}(\Omega, S)$ by $A^F(\Omega, S)$ and $X^{\infty}(\Omega, S)$ by $X^F(\Omega, S)$. For the proof, observe that if we use the abstract theory framework, it suffices by Lemma 3.15 to show that each of the classes discussed in the theorems of the section is not empty. But this is trivial, since we have shown the results for the case of $A^{\infty}(\Omega, S)$, which is a subset of $A^F(\Omega, S)$; therefore, $X^{\infty}(\Omega, S)\subset X^F(\Omega, S)$. 

\section{Topological properties of partially smooth universal Taylor series} 
In this section we will prove certain topological properties of some sets, which will allow us to deduce that there do not exist universal series with particular properties. \\ \\ 
\textbf{Definition 4.1} Let $f$ be a holomorphic function $f:\, V\rightarrow C$, where $V\subseteq C^d$ is open and $\zeta$ be a point in $\overline{V}$, such that for all differential operators of mixed partial derivatives $D$, $Df$ is continuously extendable on the point $\zeta$. We call a point $z$, $\zeta$-universal with respect to $f$, if the partial sums $S_n(f, \zeta)(z)$ are dense in $C$. We denote the set of $\zeta$-universal points with respect to $f$ by $U(f,\zeta)$. \\ \\ 
\textbf{Theorem 4.2} \textit{Let $f$ be a function $f:\,V\rightarrow C$ and $\zeta$ be a point as in Definition 4.1. Then, the set of $\zeta$-universal points is a $G_{\delta}$ subset of $C^d$.} \\ \\ 
\textit{Proof.} Let $q_n,\, n\geq 1$ be a dense sequence in $C^d$. If $G_{n,k,m}$ denotes the set $$G_{n,k,m}=\Big\{ z\,\in C^d: \vert S_k(f,\zeta)(z)-q_n\vert<\dfrac{1}{m} \Big\},$$ we easily see that $G_{n,k,m}$ is open and that 
$$U(f, \zeta)= \bigcap\limits_{n, m\geq1}\bigcup\limits_{k\geq 1}G_{n,k,m}$$ which is a $G_{\delta}$ set. \qed \\ \\ 
\textbf{Definition 4.3} Let $\Omega_i$, $1\leq i \leq d$, be simply connected domains and $\Omega=\prod\limits_{i=1}^{d}\Omega_i$. Let, also, $f:\,\Omega \rightarrow C$ be a holomorphic function on $\Omega$. We call a point $z\,\in\partial\,\Omega$ \textit{balanced}, if $Df$ can be continuously extended on the point $z$ for all differential operators of mixed partial derivatives $D$. We denote this set by $B(f)$. \\ \\ 
\textbf{Theorem 4.4} \textit{Let $f$ be as in Definition 4.3. Then the set $B(f)$ is $G_{\delta}$ in $C^d$.} \\ \\ 
\textit{Proof.} Let $D_k$, $k\geq 1$, be an enumeration of all differential operators of mixed partial derivatives. Denote by $A_{n,k}$ the set 
$$A_{n,k}=\Big\{z: \text{there exists a}\,\, \delta_z>0: \vert D_kf(s)-D_kf(t) \vert<\dfrac{1}{n}\,\, \text{whenever}\,\, s, t\in B(z, \delta_z)\cap\Omega \Big\}$$
The set $A_{n,k}$ is open and we can easily see that $B(f)=\bigcap\limits_{n,k\geq1}A_{n,k}$. \qed \\ \\ 
\textbf{Definition 4.5} Let $U$ be an open set in $C^d$, $T\subseteq \partial U$ and $W=U \cup T$. A function $f$ is \textit{$W$-bipartite universal} if 
\begin{enumerate}[label=(\roman*)]
\item $f$ is holomorphic on $\Omega$ and $Df$ is continuously entendable on $W$ for each differential operator of mixed partial derivatives $D$, 
\item for every $L$ compact subset of $W$ and any product $K$ of planar compact sets with connected complements, that is contained in $C^d-W$, and for every $h\in A_D(K)$, there exists $\lambda_n\in N$ such that $$\sup\limits_{z\,\in K,\,\zeta\in L}\big\vert S\lambda_n(f, \zeta)(z)-h(z) \big\vert \rightarrow 0 \quad \text{and} \quad \sup\limits_{z\,\in L,\,\zeta\in L}\big\vert DS\lambda_n(f, \zeta)(z)-Df(z) \big\vert \rightarrow 0$$ for every differential operator of mixed partial derivatives $D$. 
\end{enumerate} 
We denote this class by $B_u(W)$.\\ \\
\noindent
\textbf{Definition 4.6} A holomorphic function $f$ is called \textit{$W$-completely universal} if it is $W$-bipartite universal and $B(f)\cap U(f,\zeta)=\varnothing$ for every $\zeta\in W$. We denote this set by $C_u(W)$. \\ \\ 
\textbf{Remark 4.7} Theorem 3.7 implies that, when $\Omega_i$ are planar simply connected domains and $S_i$ are subsets of $\partial\,\Omega_i$, clopen in the relative topology, such that $C-(\Omega_i\cup S_i)$ is connected, then $B_u\left( \prod\limits_{i=1}^{d} (\Omega_i\cup S_i) \right)$ is a $G_{\delta}$ dense set in $X^{\infty}\left(\Omega, \prod\limits_{i=1}^{d} (\Omega_i\cup S_i)\right)$. \\ \\ 
\textbf{Corollary 4.8} \textit{There exists no $W$-bipartite universal function, where $W$ is the union of the upper half plane with the rationals of the real line.} \\ \\ 
\textit{Proof.} It is immediate by Theorem 4.4. \qed \\ \\
\textbf{Theorem 4.9} \textit{If there exists $f\in C_u(W)$, then $W\cap \partial W$ is both $F_{\sigma}$ and $G_{\delta}$.} \\ \\ 
\textit{Proof.} The proof follows from the fact that $B(f)=\left( U(f, \zeta)\right)^c$, for every $\zeta\in S$. \qed \\ \\ 
\textbf{Corollary 4.10} \textit{There exists no $S$-completely universal function $f$, where $S$ is the union of the upper half plane with the rationals of the real line, or with the irrationals of the real line.} \\ \\ 
\textit{Proof.} It is an immediate consequence of Theorem 4.9. \qed \\ \\ \\
\textit{\textbf{Acknowledgements:}} The author would like to express his gratitude towards Professor Vassili Nestoridis for the many fruitful discussions and the guidance throughout the creation of this paper.

\noindent
GIORGOS KOTSOVOLIS: Department of Mathematics, National and Kapodistrian University of Athens, Panepistimioupolis 157-84, Athens, Greece. \\ \\ 
\textit{E-mail:} \textbf{georgekotsovolis@yahoo.com}

\end{document}